\documentclass[preprint,12pt]{elsarticle}

\usepackage{amsmath,amssymb,eepic}

\renewcommand{\le}{\leqslant}
\renewcommand{\ge}{\geqslant}
\renewcommand{\sec}{\cap}
\renewcommand{\phi}{\varphi}

\makeatletter\let\ellipsearc\@arc\makeatother

\newcommand{\0}{\emptyset}
\newcommand{\1}{\mathbf{1}}
\newcommand\UNION\bigcup

\usepackage{theorem}

\theorembodyfont{\slshape}

  \newtheorem{THM}{Theorem}[section]
  \newtheorem{LEM}[THM]{Lemma}
  \newtheorem{PROP}[THM]{Proposition}
  \newtheorem{COR}[THM]{Corollary}

\theorembodyfont{\rmfamily}

  \newtheorem{DEF}[THM]{Definition}

\newif\ifQEDsign
\newcommand{\QED}{\global\QEDsigntrue\hfill$\square$}

\newenvironment{PROOF}%
    {\par\noindent\textit{Proof.}\global\QEDsignfalse}%
    {\ifQEDsign\else\QED\fi\par\bigskip\par}

\journal{Novi Sad Journal of Mathematics}

\begin{document}

\begin{frontmatter}

\title{Finite irreflexive homomorphism-homogeneous\\binary relational systems\tnoteref{THX}}
\tnotetext[THX]{Supported by the Grant No 144017 of the Ministry of Science of the republic of Serbia}

\author{Dragan Ma\v sulovi\'c, Rajko Nenadov, Nemanja \v Skori\'c}
\ead{$\{$dragan.masulovic, rajko.nenadov, nemanja.skoric$\}$@dmi.uns.ac.rs}
\address{Department of Mathematics and Informatics, University of Novi Sad\\
        Trg Dositeja Obradovi\'ca 4, 21000 Novi Sad, Serbia}

\begin{abstract}
  A structure is called homogeneous if every isomorphism between finite
  substructures of the structure extends to an automorphism of the structure.
  Recently, P.~J.~Cameron and J.~Ne\v set\v ril introduced a relaxed version
  of homogeneity: we say that a structure is homo\-mor\-phism-homogeneous if every
  homomorphism between finite substructures of the structure extends to an endomorphism
  of the structure. In this paper we characterize all finite homomorphism-homogeneous
  relational systems with one irreflexive binary relation.
\end{abstract}


\end{frontmatter}

\section{Introduction}

A structure is \emph{homogeneous} if every isomorphism between finite
substructures of the structure extends to an automorphism of the structure.
For example, finite and countably infinite homogeneous directed graphs were described in~\cite{cherlin}.
In their recent paper \cite{cameron-nesetril} the authors discuss
a generalization of homogeneity to various types of morphisms between structures,
and in particular introduce the notion of homomorphism-homogeneous structures:

\begin{DEF}[Cameron, Ne\v set\v ril \cite{cameron-nesetril}]\label{nans.def.1}
  A structure is called \emph{homo\-mor\-phism-homogeneous} if every homomorphism between
  finite substructures of the structure extends to an endomorphism of the structure.
\end{DEF}

In this short note we characterize all finite homomorphism-homogeneous
relational systems with one irreflexive binary relation.

\section{Preliminaries}
\label{nans.sec.prelim}

A binary relational system is an ordered pair $(V, E)$ where $E \subseteq V^2$ is a binary
relation on $V$. A binary relational system $(V, E)$ is \emph{reflexive} if $(x, x) \in E$ for
all $x \in V$, \emph{irreflexive} if $(x, x) \notin E$ for all $x \in V$, \emph{symmetric}
if $(x, y) \in E$ implies $(y, x) \in E$ for all $x, y \in V$ and \emph{antisymmetric}
if $(x, y) \in E$ implies $(y, x) \notin E$ for all distinct $x, y \in V$.

Binary relational systems can be thought of in terms of digraphs (hence the notation $(V, E)$).
Then $V$ is the set of \emph{vertices} and $E$ is the set of \emph{edges} of the binary relational
system/digraph $(V, E)$. Edges of the form $(x, x)$ are called \emph{loops}. If $(x, x) \in E$
we also say that \emph{$x$ has a loop}. Instead of $(x, y) \in E$ we often
write $x \to y$ and say that $x$ \emph{dominates} $y$, or that $y$ \emph{is dominated by} $x$.
By $x \sim y$ we denote that $x \to y$ or $y \to x$, while
$x \rightleftarrows y$ denotes that $x \to y$ and $y \to x$. If $x \rightleftarrows y$,
we say that $x$ and $y$ form a \emph{double edge}. We shall also say that
a vertex $x$ is \emph{incident with a double edge} if there is a vertex $y \ne x$ such that
$x \rightleftarrows y$.

Digraphs $(V, E)$ where $E$ is a symmetric binary relation on $V$ are usually referred to as \emph{graphs}.
\emph{Proper digraphs} are digraphs $(V, E)$ where $E$ is an antisymmetric binary relation.
In this paper, digraphs $(V, E)$ where $E$ is neither antisymmetric nor symmetric will be referred to
as \emph{improper digraphs}. In an improper digraph there exists a pair of distinct
vertices $x$ and $y$ such that $x \rightleftarrows y$ and another pair of distinct vertices
$u$ and $v$ such that $u \to v$ and $v \not\to u$.

A digraph $D' = (V', E')$ is a \emph{subdigraph} of a digraph $D = (V, E)$ if
$V' \subseteq V$ and $E' \subseteq E$. We write $D' \le D$ to denote that
$D'$ is isomorphic to a subdigraph of~$D$.
For $\0 \ne W \subseteq V$ by $D[W]$ we denote the digraph $(W, E \sec W^2)$
which we refer to as the \emph{subdigraph of $D$ induced by $W$}.

Vertices $x$ and $y$ are \emph{connected in $D$} if there exists a sequence of vertices
$z_1, \ldots, z_k \in V$ such that $x = z_1 \sim \ldots \sim z_k = y$. A digraph $D$ is
\emph{weakly connected} if each pair of distinct vertices of $D$ is connected in~$D$.
A digraph $D$ is \emph{disconnected} if it is not weakly connected.
A \emph{connected component of $D$} is a maximal set $S \subseteq V$ such that
$D[S]$ is weakly connected. The number of connected components of $D$ will be denoted by $\omega(D)$.

Vertices $x$ and $y$ are \emph{doubly connected in $D$} if there exists a sequence of vertices
$z_1, \ldots, z_k \in V$ such that $x = z_1 \rightleftarrows \ldots \rightleftarrows z_k = y$.
Define a binary relation $\theta(D)$ on $V(D)$ as follows:
$(x, y) \in \theta(D)$ if and only if $x = y$ or $x$ and $y$ are doubly connected.
Clearly, $\theta(D)$ is an equivalence relation on $V(D)$ and $\omega(D) \le |V(D) / \theta(D)|$.
We say that a digraph $D$ is \emph{$\theta$-connected} if $\omega(D) = |V(D) / \theta(D)|$,
and that it is \emph{$\theta$-disconnected} if $\omega(D) < |V(D) / \theta(D)|$.
Note that a $\theta$-connected digraph need not be connected, and that a $\theta$-disconnected digraph
need not be disconnected; a digraph $D$ is $\theta$-connected if every connected component of $D$ contains
precisely one $\theta(D)$-class, while it is $\theta$-disconnected if there exists a connected component of
$D$ which consists of at least two $\theta(D)$-classes. In particular, every proper digraph with at least
two vertices is $\theta$-disconnected, and every graph is $\theta$-connected.

Let $K_n$ denote the complete irreflexive graph on $n$ vertices.
Let $\1$ denote the trivial digraph with only one vertex and no edges,
and let $\1^\circ$ denote the digraph with only one vertex with a loop.
An \emph{oriented cycle with $n$ vertices} is a digraph $C_n$ whose vertices are
$1$, $2$, \ldots, $n$, $n \ge 3$, and whose edges are $1 \to 2 \to \ldots \to n \to 1$.

For digraphs $D_1 = (V_1, E_1)$ and $D_2 = (V_2, E_2)$, by $D_1 + D_2$ we denote the
\emph{disjoint union} of $D_1$ and $D_2$. We assume that $D + O = O + D = D$, where $O = (\0, \0)$
denotes the \emph{empty digraph}.
The disjoint union $\underbrace{D + \ldots + D}_k$ consisting of $k \ge 1$ copies of $D$ will be abbreviated
to $k \cdot D$. Moreover, we let $0 \cdot D = O$.

Let $D_1 = (V_1, E_1)$ and $D_2 = (V_2, E_2)$ be digraphs. We say that
$f : V_1 \to V_2$ is a \emph{homomorphism} between $D_1$ and $D_2$ and write $f : D_1 \to D_2$
if
$$
  x \to y \text{ implies } f(x) \to f(y), \text{ for all } x, y \in V_1.
$$
An \emph{endomorphsim} is a homomorphism from $D$ into itself.
A mapping $f : V_1 \to V_2$ is an \emph{isomorphism} between $D_1$ and $D_2$ if
$f$ is bijective and
$$
  x \to y \text{ if and only if } f(x) \to f(y), \text{ for all } x, y \in V_1.
$$
Digraphs $D_1$ and $D_2$ are \emph{isomorphic} if there is an isomorphism between them.
We write $D_1 \cong D_2$. An \emph{automorphsim} is an isomorphism from $D$ onto itself.

A digraph $D$ is \emph{homomorphism-homogeneous} if every homomorphism $f : W_1 \to W_2$ between
finitely induced subdigraphs of $D$ extends to an endomorphism of~$D$
(see Definition~\ref{nans.def.1}).

\section{Finite irreflexive binary relational systems}

Cameron and Ne\v set\v ril have shown in \cite{cameron-nesetril} that a finite irreflexive graph
is homomorphism-homogeneous if and only if it is isomorphic to $k \cdot K_n$ for some $k, n \ge 1$.
It was shown in \cite[Theorem 3.10]{masul-hhd} that a finite irreflexive proper digraph
is homomorphism-homogeneous if and only if it is isomorphic to
$k \cdot \1$ for some $k \ge 1$ or $k \cdot C_3$ for some $k \ge 1$.
In this section we show that these are the only finite homomorphism-homogeneous irreflexive
binary relational systems by showing that no finite irreflexive improper digraph is
homomorphism-homogeneous.

\begin{LEM}\label{nans.lem.nl-1}
  Let $D$ be a finite homomorphism-homogeneous irreflexive improper digraph.
  Then every vertex of $D$ is incident with a double edge.
\end{LEM}
\begin{PROOF}
  Let $x \rightleftarrows y$ be a double edge in $D$ and let $v$ be an arbitrary vertex of~$D$.
  The mapping
  $$
    f : \begin{pmatrix}
      x\\
      v
    \end{pmatrix}
  $$
  is a homomorphism between finitely induced subdigraphs of $D$, so it
  extends to an endomorphism $f^*$ of~$D$ by the homogeneity requirement.
  Then $x \rightleftarrows y$ implies $v = f^*(x) \rightleftarrows f^*(y)$.
\end{PROOF}

\begin{LEM}\label{nans.lem.nl-2}
  Let $D$ be a finite homomorphism-homogeneous irreflexive improper digraph and
  let $S \in V(D) / \theta(D)$ be an arbitrary equivalence class of $\theta(D)$.
  Then $D[S] \cong K_n$ for some $n \ge 2$.
\end{LEM}
\begin{PROOF}
  Lemma~\ref{nans.lem.nl-1} implies that $|S| \ge 2$ for every $S \in V(D) / \theta(D)$.

  Suppose that there is an $S \in V(D) / \theta(D)$ such that
  $D[S]$ is not a complete graph. Then there exist $u, v \in S$ such that $u \not\to v$ or $v \not\to u$.
  Let $z_1, z_2, \ldots, z_k \in V(D)$ be the shortest sequence of vertices of $D$ such that
  $$
    u = z_1 \rightleftarrows z_2 \rightleftarrows \ldots \rightleftarrows z_k = v.
  $$
  Then $k \ge 3$ since $u \not \rightleftarrows v$, and the fact that $z_1, z_2, \ldots, z_k$ is the
  shortest such sequence implies that $z_1 \not \rightleftarrows z_3$. The mapping
  $$
    f_1 : \begin{pmatrix}
      z_1 & z_3\\
      z_2 & z_3
    \end{pmatrix}
  $$
  is a homomorphism between finitely induced subdigraphs of $D$, so it
  extends to an endomorphism $f_1^*$ of~$D$ by the homogeneity requirement. Let $x_1 = f_1^*(z_2)$.
  It is easy to see that $x_1 \notin \{z_1, z_2, z_3\}$ and $x_1 \rightleftarrows y$ for all
  $y \in \{z_2, z_3\}$. Consider now the mapping
  $$
    f_2 : \begin{pmatrix}
      z_1 & z_3 & x_1 \\
      z_2 & z_3 & x_1
    \end{pmatrix}.
  $$
  which is clearly a homomorphism between finitely induced subdigraphs of $D$. It
  extends to an endomorphism $f_2^*$ of~$D$. Let $x_2 = f_2^*(z_2)$.
  Again, it is easy to see that $x_2 \notin \{z_1, z_2, z_3, x_1\}$ and that $x_2 \rightleftarrows y$ for all
  $y \in \{z_2, z_3, x_1\}$. Analogously, the mapping
  $$
    f_3 : \begin{pmatrix}
      z_1 & z_3 & x_1 & x_2 \\
      z_2 & z_3 & x_1 & x_2
    \end{pmatrix}
  $$
  is a homomorphism between finitely induced subdigraphs of $D$, so it
  extends to an endomorphism $f_3^*$ of~$D$. Let $x_3 = f_3^*(z_2)$.
  Again, $x_3 \notin \{z_1, z_2, z_3, x_1, x_2\}$ and $x_2 \rightleftarrows y$ for all
  $y \in \{z_2, z_3, x_1, x_2\}$. And so on.
  We can continue with this procedure as many times as we like, which contradicts the fact that
  $D$ is a finite digraph.
\end{PROOF}

\begin{PROP}
  There does not exist a finite homomorphism-homogeneous irreflexive improper digraph.
\end{PROP}
\begin{PROOF}
  Suppose that $D$ is a finite homomorphism-homogeneous irreflexive improper digraph.
  Then there exist vertices $x, y \in V(D)$ such that $x \to y$ and $y \not\to x$.
  Let $S = x / \theta(D)$ and $T = y / \theta(D)$. Clearly, $S \sec T = \0$.
  Let $T = \{y, t_1, \ldots, t_k\}$. Since $D[T]$ is a complete graph
  (Lemma~\ref{nans.lem.nl-2}), the mapping
  $$
    f : \begin{pmatrix}
      x & t_1 & \ldots & t_k \\
      y & t_1 & \ldots & t_k
    \end{pmatrix}
  $$
  is a homomorphism between finitely induced subdigraphs of $D$, so it
  extends to an endomorphism $f^*$ of~$D$ by the homogeneity requirement.
  Let us compute $f^*(y)$. From $f^*(t_1) \in T$ it follows that $f^*(T) \subseteq T$.
  Moreover, $f^*|_T$ is injective since there are no loops in $D$. Therefore,
  $f^*|_T : T \to T$ is a bijection. But $f^*(t_i) = t_i$ for all $i \in \{1, \ldots, k\}$,
  so it follows that $f^*(y) = y$. Now, $x \to y$ implies $f^*(x) \to f^*(y)$, that is,
  $y \to y$, which is impossible since there are no loops in~$D$.
\end{PROOF}

\begin{COR}
  Let $D$ be a finite irreflexive binary relational system. Then $D$ is homomorphism-homogeneous
  if and only if it is isomorphic to one of the following:
  \begin{enumerate}
  \item
    $k \cdot K_n$ for some $k, n \ge 1$;
  \item
    $k \cdot C_3$ for some $k \ge 1$.
  \end{enumerate}
\end{COR}

\end{document}